\documentclass[11pt]{amsart}
\usepackage{amscd,amssymb}
\usepackage[arrow,matrix]{xy}
\usepackage[colorlinks,plainpages,backref,urlcolor=blue]{hyperref}

\newtheorem{theorem}{Theorem}[section]

\newtheorem{prop}[theorem]{Proposition}

\theoremstyle{definition}
\newtheorem{remark}[theorem]{Remark}
\newtheorem*{ack}{Acknowledgments}

\newcommand{\Z}{\mathbb{Z}}

\newcommand{\C}{\mathbb{C}}

\newcommand{\T}{\mathbb{T}}

\newcommand{\E}{{\mathcal E}}

\DeclareMathOperator{\im}{im}
\DeclareMathOperator{\Hom}{Hom}
\DeclareMathOperator{\SL}{SL}

\DeclareMathOperator{\Ort}{O}

\DeclareMathOperator{\PD}{PD}

\begin{document}

\title[Which $3$-manifold groups are K\"{a}hler groups?]%
{Which $3$-manifold groups are K\"{a}hler groups?}

\author[A.~Dimca]{Alexandru Dimca}
\address{Laboratoire J.A.~Dieudonn\'{e}, UMR du CNRS 6621, 
Universit\'{e} de Nice--Sophia Antipolis, Parc Valrose,
06108 Nice Cedex 02, France}
\email{dimca@math.unice.fr}

\author[A.~Suciu]{Alexander~I.~Suciu}
\address{Department of Mathematics,
Northeastern University,
Boston, MA 02115, USA}
\email{a.suciu@neu.edu}



\begin{abstract}
The question in the title, first raised by Goldman 
and Donaldson, was partially answered by Reznikov. 
We give a complete answer, as follows: if $G$ can be realized 
as both the fundamental group of a closed $3$-manifold and 
of a compact K\"{a}hler manifold, then $G$ must be finite---%
and thus belongs to the well-known list of finite subgroups 
of $\Ort(4)$. 
\end{abstract}
\maketitle

\section{Introduction}
\label{sect:intro}

\subsection{}
As is well-known, every finitely presented group $G$ occurs as the 
fundamental group of a smooth, compact, connected, orientable 
$4$-dimensional manifold $M$.  As shown by Gompf \cite{Go}, 
the manifold $M$ can be chosen to be symplectic. Requiring a 
complex structure on $M$ is no more restrictive, as long as one 
is willing to go up to complex dimension $3$, see Taubes \cite{Ta}. 

Suppose now $G$ is the fundamental group of a compact 
K\"{a}hler manifold $M$.  Groups arising this way are called 
{\em K\"{a}hler groups} (or, {\em projective groups}, if $M$ is 
actually a smooth projective variety).  The K\"{a}hler 
condition puts strong restrictions on what $G$ can be.  
For instance, the first Betti number, $b_1(G)$, must 
be even, by classical Hodge theory. Moreover, $G$ must 
be $1$-formal, by work of Deligne, Griffiths, Morgan, and 
Sullivan \cite{DGMS}. Also, $G$ cannot split non-trivially 
as a free product, by a result of Gromov \cite{Gr}.  
On the other hand, every finite group is a projective  group, 
by a classical result of Serre \cite{Se}. We refer to \cite{ABCKT} 
for a comprehensive survey of K\"{a}hler groups, and 
to the recent work of Delzant--Gromov \cite{DG}, 
Napier--Ramachandran \cite{NR}, and Delzant \cite{De}  
for further geometric restrictions imposed by the K\"{a}hler 
condition on a group $G$.  

Requiring that $M$ be a $3$-dimensional compact, connected 
manifold also puts severe restrictions on $G=\pi_1(M)$.  
For example, if $G$ is abelian, then $G$ is either $\Z/n\Z$, 
$\Z$, $\Z\oplus \Z_2$, or $\Z^3$, see \cite{Ja}. 

\subsection{}
A natural question---raised by Goldman and Donaldson 
in 1989, and independently by Reznikov in 1993---is then:  
what are the $3$-manifold groups which are K\"{a}hler groups? 

In \cite{Re}, Reznikov proved the following result, which 
Simpson \cite{Si04} calls ``one of the deepest restrictions" 
on the homotopy types that may occur for K\"{a}hler manifolds:
{\em Let $M$ be an irreducible, atoroidal $3$-manifold, and suppose 
there is a homomorphism $\rho\colon \pi_1(M)\to \SL(2,\C)$ 
with Zariski dense image. Then $G=\pi_1(M)$ is not a 
K\"{a}hler group.}  The same conclusion was reached by 
Hern\'{a}ndez-Lamoneda in \cite{HL}, under the assumption 
that $M$ is a geometrizable $3$-manifold, with all pieces 
hyperbolic.   

In this note, we answer the above question for all $3$-manifold 
groups, as follows.  

\begin{theorem}
\label{thm:main}
Let $G$ be the fundamental group of a compact, connected 
$3$-manifold. If $G$ is a K\"{a}hler group, then $G$ is finite.   
\end{theorem}

By the  $3$-dimensional spherical space-form conjecture, now 
established by Perelman \cite{Pe2, Pe3}, a closed 
$3$-manifold $M$ has finite fundamental group if and 
only if it admits a metric of constant positive curvature 
(for a detailed proof, see Morgan and Tian \cite[Corollary 0.2]{MT}). 
Thus, $M=S^3/G$, where $G$ is a finite subgroup of $\Ort(4)$, 
acting freely on $S^3$.  The list of such finite groups (essentially  
due to Hopf) is given by Milnor in \cite{Mi}. 

\subsection{}
The paper is organized as follows.  In \S\ref{sect:cjl}, 
we discuss the characteristic and resonance varieties 
of a group $G$, and two notions of isotropicity. 
In \S\ref{sect:ist}, we recall the Isotropic Subspace 
Theorem of Catanese, and a correspondence due to 
Beauville. In \S\ref{sect:resk}, we use these tools to 
prove a key result, tying the first resonance variety of 
a K\"{a}hler manifold to the rank of the cup-product 
map in low degrees. In \S\ref{sect:res3}, we investigate 
the first resonance variety of a closed, oriented 
$3$-manifold; Poincar\'{e} duality and properties 
of Pfaffians yield a very different conclusion in this setting.  

All this works quite well, provided the first Betti number 
of $G$ is positive. To deal with the remaining case, we need 
two theorems of Reznikov and Fujiwara, relating the K\"{a}hler, 
respectively the $3$-manifold condition on a group to Kazhdan's 
property $T$; we recall those in \S\ref{sect:h2k}.  Finally, 
we put everything together in \S\ref{sect:k3}, and give a proof 
of Theorem~\ref{thm:main}.

A natural question arises out of this work:  Which $3$-manifold 
groups are quasi-K\"{a}hler?  (A group $G$ is {\em quasi-K\"{a}hler}\/ 
if $G=\pi_1(M\setminus D)$, where $M$ is a compact K\"{a}hler 
manifold and $D$ is a divisor with normal crossings.) 
We have some partial results in this direction; 
those results will be presented elsewhere. 

\section{Cohomology jumping loci and isotropic subspaces}
\label{sect:cjl}

\subsection{}
Let $X$ be a connected CW-complex with finitely many cells 
in each dimension. Let $G=\pi_1(X)$ be the fundamental group 
of $X$, and $\T =\Hom (G, \C^*)$ its character variety. 
Every character $\rho\in \T$ determines a rank $1$ local 
system, $\C_{\rho}$, on $X$.   The {\em characteristic varieties}\/ 
of $X$ are the jumping loci for cohomology with coefficients in 
such local systems:
\begin{equation}
\label{eq:cv}
V^i_d(X)=\{\rho \in \T
\mid \dim H^i(X,\C_{\rho})\ge d\}. 
\end{equation}

The varieties $V_d(X)=V^1_d(X)$ depend only on $G=\pi_1(X)$, 
so we sometimes denote them as $V_d(G)$.  

\subsection{}
Consider now the cohomology algebra $A=H^* (X,\C)$.  
Left-mul\-tiplication by an element $x\in A^1$ yields a 
cochain complex $(A, x)\colon A^0\xrightarrow{x} 
A^1\xrightarrow{x} A^2 \to \cdots$.  
The {\em resonance varieties}\/ of $X$ are the jumping 
loci for the homology of this complex:
\begin{equation}
\label{eq:rv}
R^i_d(X)=\{x \in A^1 \mid 
\dim H^i(A,x) \ge  d\}.
\end{equation}

The varieties $R_d(X)=R^1_d(X)$ depend only on 
$G=\pi_1(X)$, so we sometimes denote them by $R_d(G)$. 
By definition, an element $x\in A^1$ belongs to $R_d(X)$ 
if and only if there exists a subspace $W\subset A^1$ of 
dimension $d+1$ such that $x\cup y=0$, for all $y\in W$.  

Fix bases $\{e_1,\dots ,e_n\}$ for $A^1$ and $\{f_1,\dots ,f_m\}$ 
for $A^2$.  Writing the cup-product as 
$e_i\cup e_j = \sum_{k=1}^m \mu_{i,j,k} f_k$, 
we may define an $m\times n$ matrix $\Delta$ of linear 
forms in variables $x_1,\dots , x_n$, with entries  
\begin{equation}
\label{eq:delta}
\Delta_{k,j} =  \sum_{i=1}^n \mu_{i,j,k} x_i.  
\end{equation}
It is readily seen that $R_d(X)=V(E_{d} (\Delta))$, 
where $E_{d}$ denotes the ideal of $(n-d)\times (n-d)$ 
minors.  Note also that $x\cup x=0$, for all $x\in A^1$ implies 
$\Delta\cdot \vec{x}=0$, where $\vec{x}$ is the column 
vector with entries $x_1,\dots ,x_n$. 

\subsection{}
Foundational results on the structure of the cohomology support 
loci for local systems on compact K\"{a}hler manifolds 
were obtained by Beauville \cite{Be}, Green--Lazarsfeld 
\cite{GL}, Simpson \cite{Si92}, and Campana \cite{Cm01}: 
if $G$ is the fundamental group of such a manifold, then 
$V_d(G)$ is a union of (possibly translated) subtori of the 
algebraic group $\T$. 

In addition, Theorem~A from \cite{DPS05} 
establishes a strong relationship between the characteristic 
and resonance varieties of a K\"{a}hler group $G$:  the 
tangent cone to $V_d(G)$ at the identity of $\T$ 
equals $R_d(G)$, for all $d\ge 1$. 

\subsection{}
A non-zero subspace $E\subset H^1(X,\C)$ is {\em (totally) isotropic}\/ 
if the restriction of the cup-product map $\cup_X\colon H^1(X,\C) \wedge 
H^1(X,\C) \to H^2(X,\C)$ to $E\wedge E$ is identically zero.  By analogy, 
we say $E$ is {\em $1$-isotropic}\/ if the restriction of $\cup_X$ to 
$E\wedge E$ has $1$-dimensional image.

Note that the these properties of $E$ depend only on $G=\pi_1(X)$. 
Indeed, let $h\colon X\to K(G,1)$ be a classifying map.  
Then $h_*\colon H_1(X,\Z)\to H_1(G,\Z)$ is an isomorphism, 
and $h_*\colon H_2(X,\Z)\to H_2(G,\Z)$ is an epimorphism. 
Using Kronecker duality and the functoriality of the cup-product, 
it is readily seen that $E$ is an ($1$-) isotropic subspace of 
$H^1(G,\C)$ for $\cup_G$
if and only if $h^*(E)$ is an ($1$-) 
isotropic subspace of $H^1(X,\C)$ for $\cup_X$. 

\section{The Isotropic Subspace Theorem}
\label{sect:ist}

By a {\em fibration}\/ we mean a surjective morphism $f\colon M \to N$ 
with connected fibers between two compact complex manifolds $M$ and $N$. 
Two fibrations $f\colon M \to C$ and $f'\colon M\to C'$ over projective curves 
$C$ and $C'$ are said to be {\em equivalent}\/ if there is an isomorphism 
$\phi\colon C \to C'$ such that $f'=\phi \circ f$. We denote by $\E (M)$ the 
set of equivalence classes of all these fibrations.

Beauville's work \cite{Be} establishes a bijection 
between the set $\E (M)$ and the set of irreducible components 
of the first characteristic variety $V_1(M)$ passing through the identity 
of the algebraic group $\T=\Hom(\pi_1(M),\C^*)$. In particular, 
the set $\E (M)$ must be finite.

The Isotropic Subspace Theorem, due to Catanese \cite[Theorem 1.10]{C1}, 
establishes a relation between the set of equivalence classes of fibrations 
of a K\"{a}hler manifold $M$ over curves of genus $g\ge 2$, and  
the maximal isotropic subspaces in $H^1(M,\C)$. 

\begin{theorem}[Catanese \cite{C1}]
\label{thm:IST}
Let $M$ be a compact K\"{a}hler manifold. Then, for any maximal isotropic 
subspace $E \subset H^1(M,\C)$ of dimension $g\ge 2$, there is a fibration 
$f\colon M \to C$ onto a smooth curve of genus $g$ and a maximal isotropic 
subspace $E' \subset H^1(C,\C)$ such that $E=f^*E'$.
\end{theorem}

For more information on this correspondence, see \cite{C2}. 

\section{The first resonance variety of a K\"{a}hler manifold}
\label{sect:resk}

\begin{theorem}
\label{keythm} 
Let $M$ be a compact K\"{a}hler manifold with $b_1(M)\ne 0$.  
If $R_1(M)= H^1(M,\C)$, then $H^1(M,\C)$ is $1$-isotropic.  
\end{theorem}

\begin{proof}
By Hodge theory, we must have $b_1(M) \ge 2$. 
The equality $R_1(M)= H^1(M,\C)$ says that, for any 
non-zero cohomology class $x \in H^1(M,\C)$, there is a class 
$y \in H^1(M,\C)\setminus \C \cdot x$ such that $x \cup y=0$. 
Consequently, the vector space spanned by $x$ and $y$ is a  
($2$-dimensional) isotropic subspace containing $x$. 

Let $U_x$ be a maximal isotropic subspace of $H^1(M,\C)$ 
containing $x$; we must then have $\dim U_x \ge 2$.  
Thus, by  Theorem \ref{thm:IST}, there is a fibration 
$f_x\colon M \to C_x$ onto a smooth projective curve $C_x$ 
of genus $g_x=\dim U_x$, with $x \in f_x^*(H^1(C_x,\C))$. 

Recall now that the set $\E(M)$ of equivalence classes of fibrations  
of $M$ over curves of genus at least $2$ is finite. Thus, we may write 
the first cohomology group of $M$ as a finite union of linear 
subspaces,
\begin{equation}
\label{eq:bigcup}
H^1(M,\C)=\bigcup_{[f] \in \E(M)}f^*(H^1(C_f,\C)), 
\end{equation}
where $f=f_x$ for some $x\in H^1(M,\C)$, and $C_{f}:=C_x$. 
This is possible only if there is a fibration $f_1\colon M \to C_1$ 
such that $H^1(M,\C)=f_1^*(H^1(C_1,\C))$.

Since $f_1$ is a fibration, the induced morphism 
$f_1^*\colon H^1(C_1,\C) \to H^1(M,\C)$ 
is injective. The defining property of $f_1$ implies that 
$f_1^*\colon H^1(C_1,\C) \to H^1(M,\C)$ is an isomorphism.

On the other hand, the induced morphism 
$f_1^*\colon H^2(C_1,\C) \to H^2(M,\C)$ is also injective. 
To prove this claim, first note that any cohomology 
class in $H^1(M,\C)$ is primitive. Using the Hodge-Riemann 
bilinear relations, see e.g.~\cite[p.~123]{GH}, it follows that, 
for any non-zero $(1,0)$-class $\alpha \in H^1(M,\C)$, the 
product $\beta= \sqrt{ -1}\, \alpha \cup \overline \alpha$ 
is a non-zero, real, $(1,1)$-class in $H^2(M,\C)$. 
Since $f_1^*\colon H^1(C_1,\C)\to H^1(M,\C)$ is an isomorphism, 
there is an element $a \in H^1(C_1,\C)$ such that $f_1^*(a)=\alpha$. 
Hence, $f_1^*(\sqrt{ -1}\, a \wedge \overline a)=\beta$, and the  
claim is proved.

Consider now the commuting diagram
\begin{equation}
\label{eq:cup}
\xymatrix{
H^1(M,\C)\wedge  H^1(M,\C) \ar^(.63){\cup_M}[rr]&&  H^2(M,\C) \\
H^1(C_1,\C)\wedge  H^1(C_1,\C) \ar_{f_1^*\wedge f_1^*}[u] 
\ar^(.63){\cup_{C_1}}[rr]&& H^2(C_1,\C)\ar_{f_1^*}[u]
}
\end{equation}
As we saw above, the left arrow is an isomorphism, 
and the right one is an injection.  
Since $\cup_{C_1}$ surjects onto $H^2(C_1,\C)=\C$, we conclude 
that $\cup_M$ has $1$-dimensional image. 
\end{proof}

\begin{remark}
\label{rem:DPS}
An alternate way to prove Theorem \ref{keythm} is by using 
the much more general Theorem~B from \cite{DPS05}, which 
guarantees that {\em every}\/ positive-dimensional component 
of $R_1(M)$ is an $1$-isotropic subspace of $H^1(M, \C)$. 
This is the argument we had in an earlier version of this paper; 
at the urging of one of the referees, we came up with the above, 
more self-contained proof. 
\end{remark}

\section{The first resonance variety of a $3$-manifold}
\label{sect:res3}

Let $M$ be a compact, connected, orientable $3$-manifold.  
Fix an orientation on $M$, that is, pick a generator 
$[M]\in H^3(M,\Z)\cong\Z$. With this choice, the cup 
product on $M$ determines an alternating $3$-form 
$\mu=\mu_M$ on $H^1(M,\Z)$, given by 
\begin{equation}
\label{eq:mu}
\mu(x,y,z) = \langle x\cup y\cup z , [M]\rangle,
\end{equation}
where $\langle\, , \rangle$ is the Kronecker pairing. 
In turn, the cup-product map $\cup_M\colon H^1(M,\Z) \wedge 
H^1(M,\Z) \to H^2(M,\Z)$ is determined by $\mu$, via 
$\langle x\cup y , \gamma \rangle = \mu (x,y,z)$, 
where $z=\PD(\gamma)$ is the Poincar\'{e} dual of 
$\gamma \in H_2(M,\Z)$. 

Now fix a basis $\{e_1,\dots ,e_n\}$ for $H^1(M,\C)$, 
and choose as basis for $H^2(M,\C)$ the set 
$\{e^{\vee}_1,\dots ,e^{\vee}_n\}$, where 
$e^{\vee}_i$ denotes the Kronecker dual 
of the Poincar\'e dual of $e_i$.  Then 
\begin{equation}
\label{eq:muform}
\mu(e_i,e_j,e_k )= \langle \sum_{1\le m \le n} \mu_{i,j,m}e^{\vee}_m, 
\PD(e_k)\rangle= \mu_{i,j,k}.
\end{equation}
Recall from \eqref{eq:delta} the $n\times n$ matrix 
with entries $\Delta_{k,j} =  \sum_{i=1}^n \mu_{i,j,k} x_i$.  
Since $\mu$ is an alternating form, $\Delta$ 
is a skew-symmetric matrix.

\begin{prop}
\label{prop:res3}
Let $M$ be a closed, orientable $3$-manifold.  
Then:
\begin{enumerate}
\item  \label{r1}
$H^1(M,\C)$ is not $1$-isotropic. 
\item \label{r2}
If $b_1(M)$ is even, then $R_1(M)=H^1(M,\C)$.
\end{enumerate}
\end{prop}

\begin{proof}
To prove \eqref{r1}, suppose $\dim \im (\cup_M)=1$.  
This means there is a hyperplane $E\subset H:=H^1(M,\C)$ 
such that $x\cup y\cup z= 0$, for all $x,y \in H$ and $z\in E$.  
Hence, the skew $3$-form $\mu\colon \bigwedge^3 H \to \C$ 
factors through a skew $3$-form $\bar\mu\colon 
\bigwedge^3 (H/E) \to \C$. But $\dim H/E=1$ forces 
$\bar\mu=0$, and so $\mu=0$, a contradiction. 

To prove \eqref{r2},  recall $R_1(M)=V(E_1(\Delta))$. 
Since $\Delta$ is a skew-symmetric matrix of even size, 
it follows from Buchsbaum--Eisenbud \cite[Corollary 2.6]{BE} 
that $V(E_1(\Delta))=V(E_0(\Delta))$, 
see \cite[eq.~(6.9)]{CS06}. But  $\Delta\cdot \vec{x}=0$ 
implies $\det \Delta =0$, and so $V(E_0(\Delta))=H$. 
\end{proof}

\begin{remark}
\label{rem:bp}
As  noted by S.~Papadima, the following holds.  
Suppose $M$ is a closed, orientable $3$-manifold, with 
$b_1(M)$ odd.  Then, $R_1(M)\ne H^1(M,\C)$ if and only 
if $\mu_M$ is generic, in the sense of \cite{BP}. 
\end{remark}

\section{Kazhdan's property $T$}
\label{sect:h2k}

The following question is due to J.~Carlson and D.~Toledo 
(see J.~Koll\'{a}r \cite{Kl}):  For a K\"{a}hler group $G$, is 
$b_2(G)\ne 0$? This question was answered in the affirmative 
by A.~Reznikov in \cite{Re}, under an additional assumption, 
as follows.

\begin{theorem}[Reznikov \cite{Re}]
\label{thm:rez}
Let $G$ be a K\"{a}hler group. If $G$ does not satisfy 
Kazhdan's property $T$, then $b_2(G)\ne 0$. 
\end{theorem}

Recall that a discrete group $G$ satisfies Kazhdan's property 
$T$ (for short, $G$ is a Kazhdan group) if and only if 
$H^1(G, \mathcal{H})= 0$, for all orthogonal or unitary 
representations of $G$ on a Hilbert space $\mathcal{H}$, 
see de~la~Harpe and Valette \cite[p.~47]{HV}.  In particular, 
if $b_1(G)\ne 0$, then $G$ is not Kazhdan. (For a simple 
proof of Theorem \ref{thm:rez} in this case, see \cite{JR}.) 

We will also need the following relationship between 
$3$-manifold groups and Kazhdan's property $T$, 
established by K.~Fujiwara in \cite{Fu}.

\begin{theorem}[Fujiwara \cite{Fu}]
\label{thm:fuji}
Let $G$ be the fundamental group of a closed, orientable 
$3$-manifold. If $G$ satisfies Kazhdan's property $T$, 
then $G$ is finite. 
\end{theorem}

In fact, the theorem is valid for any subgroup $G< \pi_1(M)$, 
where $M$ is a compact (not necessarily boundaryless), connected, 
orientable $3$-manifold.   Fujiwara further assumes that each piece 
of the canonical decomposition of $M$ along embedded spheres, 
disks and tori admits one of the eight geometric structures in the 
sense of Thurston, but this is now guaranteed by the work of Perelman 
\cite{Pe2, Pe3}.

\section{K\"{a}hler $3$-manifold groups}
\label{sect:k3}

We are now in position to prove Theorem \ref{thm:main} 
from the Introduction.  

Let $G$ be the fundamental group of a compact, connected 
$3$-manifold $M$.  Suppose $G$ is a K\"{a}hler group, 
and $G$ is not finite.  

{\em Step 1.} A finite-index subgroup of a K\"{a}hler group is 
again a K\"{a}hler group (see \cite[Example 1.10]{ABCKT}). 
Passing to the orientation double cover of $M$ if 
necessary, we may as well assume $M$ is orientable. 

{\em Step 2.} 
Since $G$ is an infinite, orientable $3$-manifold group, $G$ is 
not Kazhdan, by Fujiwara's Theorem \ref{thm:fuji}.  Since $G$ 
is K\"{a}hler and not Kazhdan, $b_2(G)\ne 0$, by Reznikov's 
Theorem \ref{thm:rez}. 

{\em Step 3.} 
Since $b_2(M) \ge b_2(G)$, we must also have  
$b_2 (M) \ne 0$.  By Poincar\'{e} duality, $b_1(M) = b_2 (M)$. 
Hence, $b_1(G)=b_1(M)$ is not zero.  

{\em Step 4.} 
Since $G$ is K\"{a}hler, $b_1(G)$ must be even. 
Since $M$ is a closed, orientable $3$-manifold with $G=\pi_1(M)$, 
Proposition \ref{prop:res3} tells us that $R_1(G)=H^1(G,\C)$ 
and $H^1(G,\C)$ is not $1$-isotropic.  Since, on the other 
hand, $G$ is K\"{a}hler, Theorem \ref{keythm} tells us 
that $b_1(G)=0$. 

Our assumptions have led us to a contradiction.  Thus, the 
Theorem is proved. 

\begin{ack}
This work was done during the second author's visit at 
Universit\'{e} de Nice--Sophia Antipolis in September, 2007. 
He thanks the Laboratoire Jean A. Dieudonn\'{e}
for its support and hospitality during his stay in Nice, France.
Both authors thank the referees for constructive comments 
and helpful suggestions that led to improvements in the 
presentation and content of the paper.
\end{ack}

\vspace*{-3pc}

\newcommand{\arxiv}[1]
{\texttt{\href{http://arxiv.org/abs/#1}{arXiv:#1}}}

\renewcommand{\MR}[1]
{\href{http://www.ams.org/mathscinet-getitem?mr=#1}{MR#1}}

\end{document}